\documentclass{article}
\usepackage{amsmath,times}
\usepackage{amsthm}
\usepackage{amssymb}
\usepackage{latexsym}
\usepackage{amscd}
\usepackage[english]{babel}
\usepackage[latin1]{inputenc}
\newtheorem{theorem}{Theorem}[section]
\newtheorem{lemma}[theorem]{Lemma}

\newtheorem{definition}[theorem]{Definition}

\begin{document} 

\title{Convergence of Bergman measures for high powers of a line bundle}
\author{ROBERT BERMAN, DAVID WITT NYSTRÖM}
\maketitle

\begin{abstract}
Let $L$ be a holomorphic line bundle on a compact complex manifold $X$ of dimension $n,$ and let $e^{-\varphi}$ be a continuous metric on $L.$ Fixing a measure $d\mu$ on $X$ gives a sequence of Hilbert spaces consisting of holomorphic sections of tensor powers of $L.$ We prove that the corresponding sequence of scaled Bergman measures converges, in the high tensor power limit, to the equilibrium measure of the pair $(K,\varphi),$ where $K$ is the support of $d\mu,$ as long as $d\mu$ is stably Bernstein-Markov with respect to $(K,\varphi).$ Here the Bergman measure denotes $d\mu$ times the restriction to the diagonal of the pointwise norm of the corresponding orthogonal projection operator. In particular, an extension to higher dimensions is obtained of results concerning random matrices and classical orthogonal polynomials.
\end{abstract}

\section{Introduction}

Let $L$ be a holomorphic line bundle on a compact complex manifold $X.$ A hermitian metricon $L,$ locally represented as $h=e^{-\varphi}$ where $\varphi$ will be called the weight of the metric, together with a positive Borel measure $d\mu$ on $X,$ gives a scalar product on the space $H^0(X,L)$ of holomorphic sections of $L,$  by letting $$||s||_{\varphi}^2:=\int_X |s|^2e^{-\varphi}d\mu,$$ which we will assume to be nondegenerate. One can choose an orthonormal basis $\{s_i\}$ for $H^0(X,L)$ with respect to this scalar product, and define the Bergman function as $$B_{\varphi}(z):=\sum_i |s_i(z)|^2e^{-\varphi}.$$ This does not depend on the particular choice of orthonormal basis, and in fact the Bergman function could be defined invariantly as the pointwise norm on the diagonal of the Bergman kernel $K_{\varphi}(z,w),$ where $K_{\varphi}$ represents the orthogonal projection $C^{\infty}(X,L)\rightarrow H^0(X,L)$.

The Bergman function times the measure $d\mu$ gives a measure, $$B_{\varphi}(z)d\mu(z),$$ called the Bergman measure. Taking tensor powers of the line bundle $L,$ which we will write additively as $kL,$ with the corresponding multiplied weight $k\varphi,$ one thus gets a sequence of measures on $X,$ $$B_{k\varphi}(z)d\mu(z).$$ If scaled by a factor $1/k^n,$ this sequence is weakly bounded, and one is interested in the possible weak convergence of the scaled Bergman measures, $$\frac{B_{k\varphi}(z)}{k^n}d\mu(z).$$

In this paper we will prove that the limit does exist and is given by the equilibrium measure associated to $(K,\varphi),$ where $K$ denotes the support of $d\mu,$ as long as $d\mu$ is stably Bernstein-Markov with respect to $(K,\varphi).$ This essentially proves the conjectures made in \cite{Bloom2} and \cite{Berman2}. The proof relies heavily on the very recent work \cite{Berman}. Before giving the precise statement of the theorem to be proved we will review previous results.

In the case when $\varphi$ is smooth with stricly positive curvature form $dd^c\varphi$ (and hence $L$ is ample), and $d\mu$ the measure induced by a smooth volume form $\omega_n,$ Tian (see e.g. \cite{Zelditch2}) and Bouche (see \cite{Bouche}) independently showed that the scaled Bergman measures $B_{k\varphi}\omega_n /k^n$ converge weakly to the Monge-Ampère measure of $\varphi,$ $$\textrm{MA}(\varphi):=\frac{(dd^c \varphi)^n}{n!}.$$

When $L$ is an arbitrary line bundle and $\varphi$ any smooth weight on $L$ it was recently shown in \cite{Berman3} that the following weak convergence holds: $$\frac{B_{k\varphi}}{k^n}\omega_n \stackrel{*}{\longrightarrow} \textrm{MA}(P(\varphi)),$$ where
\begin{equation} \label{Equation10}
P(\varphi):=\sup\{\psi, \psi \textrm{ psh weight on } L, \psi \leq \varphi \textrm{ on X}\}.
\end{equation}
Here a plurisubharmonic weight denotes a weight with positive curvature current. The Monge-Ampère operator $\textrm{MA}$ is defined according to the classical work of Bedford-Taylor on any Zariski open set where $P(\varphi)$ is locally bounded and then extended by zero to all of $X$ (cf. \cite{Berman3}).

Next, consider the general situation where the weight is merely continuous and $d\mu$ is a positive measure, whose support will be denoted by $K.$ 

A classical situation is obtained by setting $X=\mathbb{P}^n$, $L=\mathcal{O}(1),$ with $K$ compactly included in the affine piece $\mathbb{C}^n.$ An element in $H^0(kL)$ may then be represented by a complex polynomial $p_k$ of total degree at most $k$ and the corresponding norm as $$\int_K |p_k|^2e^{-k\varphi}d\mu$$. One may for example let $K$ be a smooth domain, its boundary, or a smooth subset in the totally real set $\mathbb{R}^n,$ and $d\mu$ the one induced by the Lesbegue measure on $\mathbb{C}^n.$ There is also an unbounded variant of this setting where $K$ is allowed to be an unbounded set in $\mathbb{C}^n,$ but where the weight function $\varphi$ is required to satisfy the growth assumption $$\varphi(z)\geq (1+ \varepsilon)\log |z|^2,$$ for some $\varepsilon>0.$ \footnote{Hence $\varphi$ does not extend to a continuous metric on $\mathcal{O}(1) \rightarrow \mathbb{P}^n$.} This situation appears, for example, naturally in random matrix theory, where the Bergman measure $B_{k\varphi}d\mu$ represents the expected distribution of eigenvalues of a random matrix (see \cite{Berman4}, \cite{Deift}, \cite{Hedenmalm} and \cite{Saff}) of rank $k+1.$ Note that the Bergman function $B_{k\varphi}(z)$ often is refered to as the Christoffel-Darboux function in the classical literature about orthogonal polynomials. 

In the special case when $X=\mathbb{P}^1,$ $L=\mathcal{O}(1),$ and $K$ a compact subset of the affine piece $\mathbb{C}^1$ it was shown very recently by Bloom-Levenberg that if the measure $d\mu$ has the Bernstein-Markov property, i.e. for every $\delta>0$ there exist a constant $C$ such that $$\sup_{z\in K}\{|s(z)|^2e^{-k\varphi}\}\leq Ce^{\delta k}||s||_{L^2(e^{-k\varphi}d\mu)}^2,$$ for all $k$ and all $s\in H^0(kL),$ then $$\frac{B_{k\varphi}}{k^n}d\mu \stackrel{*}{\longrightarrow} \textrm{MA}(P_K(\varphi)).$$ Here $P_K(\varphi)$ is defined as the upper semicontinuous (usc) regularization of the weight obtained by replacing $X$ in formula (\ref{Equation10}) by the compact set $K,$ assumed to be locally nonpluripolar. The corresponding measure $\textrm{MA}(P_K(\varphi))$ is called the equilibrium measure of $(K,\varphi).$ See also \cite{Berman2} for results concerning the case when $K$ is a pseudoconcave domain, and \cite{Zelditch} concerning analytic domains in $\mathbb{C}^1,$ where the latter reference is based on classical results of Szegö \cite{Szego} and Carleman \cite{Carleman} concerning orthogonal polynomials.

It was conjectured by Bloom-Levenberg in (\cite{Bloom2}) that the weak convergence holds also in $\mathbb{P}^n$ for $n>1,$ under the Bernstein-Markov assumption on $d\mu.$ The corresponding conjecture in the more general line bundle setting was made in \cite{Berman2}.

In this paper we prove this under slightly stronger assumptions, namely that the measure $d\mu$ is stably Bernstein-Markov with respect to $(K, \varphi),$ meaning that not only is $d\mu$ Bernstein-Markov with respect to $(K,\varphi),$ but also with respect to $(K, \varphi+\varepsilon u)$ for any small perturbation  $\varphi+\varepsilon u$ of the weight $\varphi.$

\begin{theorem} \label{Theorem1}
Let $X$ be a compact complex manifold, $L$ a line bundle over $X$ with a continuous weight $\varphi,$ $K$ a compact subset of $X,$ and $d\mu$ a measure stably Bernstein-Markov with respect to $(K, \varphi)$. Then the sequence of Bergman measures $B_{k\varphi}d\mu/k^n$ converges weakly to the equilibrium measure $\textrm{MA}(P_K(\varphi)).$
\end{theorem}

The corresponding result also holds in the unbounded setting in $\mathbb{C}^n$ refered to above with essentially the same proof (compare remark 9.2 in \cite{Berman}).

In all the classical cases refered to above, the measure $d\mu$ will be stably Bernstein-Markov, see \cite{Berman}, \cite{Bloom2} and references therein.

As pointed out above, the proof relies heavily on the very recent work \cite{Berman}. The starting point is the fact that the Bergman measure represents the differential of a certain functional $F_k$ defined on the affine space of all continuos weights (for a fixed compact set $K$ in $X$). It was shown in \cite{Berman} that $f_k$ converges to a concave functional $F$ with continuous Frechet differential, under the assumption that $d\mu$ has the Bernstein-Markov property. Moreover, the differential of $F$ was shown to be represented by the equilibrium measure $\textrm{MA}(P_K(\varphi)).$ The new key observation that we make is that $F_k$ is in fact concave for any $k$. Then an elementary calculus lemma gives the convergence on the level of derivatives and hence concludes the proof of the theorem.

It is a pleasure to thank Bo Berndtsson and Sébastien Boucksom for fruitful discussions related to the topic of the present paper. The authors are also grateful to the Mittag-Leffler institute (Stockholm) and Institut Fourier (Grenoble) where parts of this work was carried out.

\section{Proof of Theorem 1.1}

In (\cite{Berman}) Berman-Boucksom introduced the notion of a relative capacity $\mathcal{E}_0(P_K(\varphi))$ of a weighted compact set $(K,\varphi)$. This capacity has the property that if $u$ is a smooth function, and one differentiates the perturbed capacity $\mathcal{E}_0(P_K(\varphi+\varepsilon u))$  with respect to $\varepsilon,$ one gets that $$\mathcal{E}_0(P_K(\varphi+\varepsilon u))_{\varepsilon}'=\int_M u \textrm{MA}(P_K(\varphi+\varepsilon u)).$$ The main result in (\cite{Berman}) states that if $d\mu$ is Bernstein-Markov with respect to $(K,\varphi),$ then the relative capacity is given as a limit of logarithmic volumes, $$\mathcal{E}_0(P_K(\varphi))=\lim_{k \to \infty}\frac{(n+1)!}{2k^{n+1}}\log \textrm{vol } \mathcal{B}^2(d\mu, k\varphi).$$ Here $\textrm{vol } \mathcal{B}^2(d\mu, k\varphi)$ means the volume of the unit ball in $H^0(kL)$ with respect to the $L^2(e^{-k\varphi}d\mu)$-norm, where the volume is computed relative the unique Haar measure on $H^0(kL)$ which gives volume one to the unit ball determined by a fixed but arbitrary reference weight $\varphi_0$. This in turn establishes a connection to the Bergman function. If one differentiates the function $\log \textrm{vol }\mathcal{B}^2(d\mu, \varphi + \varepsilon u),$ it is not hard to see (\cite{Berman}) that $$(\log \textrm{vol }\mathcal{B}^2(d\mu, \varphi + \varepsilon u))_{\varepsilon}'=\int_M u B_{\varphi+\varepsilon u}d\mu.$$

\begin{definition}
We say that $u$ is a Bernstein-Markov direction to $(K, \varphi, d\mu)$ if for all small $\varepsilon,$ $d\mu$ is Bernstein-Markov with respect to $(K, \varphi+\varepsilon u)$, as it was defined in the introduction.
\end{definition}

\begin{lemma} \label{Lemma2}
Let $u$ be a smooth function on $X$, and define the function $f$ by $$f(\varepsilon):= \log \textrm{vol } \mathcal{B}^2(d\mu, \varphi+\varepsilon u).$$ Then $f'' \leq 0$.
\end{lemma}

\begin{proof}
We know that $$f'(\varepsilon)=\int u B_{\varphi + \varepsilon u}d\mu,$$ where $B_{\varphi + \varepsilon u}$ denotes the Bergman function.

Let $K_{\varphi}(z,w)$ be the Bergman kernel. Since the kernel is holomorphic in the first variable, the reproducing property gives that 
\begin{equation} \label{Equation3}
K_{\varphi}(z,z)=\int K_{\varphi}(w,z)K_{\varphi}(z,w)e^{-\varphi(w)}d\mu(w).
\end{equation}
If we assume that the weight $\varphi$ depends on a parameter $\varepsilon$, differentiating with respect to $\varepsilon$ gives 
\begin{eqnarray} \label{Equation2}
K_{\varphi}(z,z)_{\varepsilon}'=\int K_{\varphi}(w,z)_{\varepsilon}'K_{\varphi}(z,w)e^{-\varphi(w)}d\mu(w)+\nonumber \\+\int K_{\varphi}(w,z)K_{\varphi}(z,w)_{\varepsilon}'e^{-\varphi(w)}d\mu(w)-\nonumber \\-\int \varphi(w)_{\varepsilon}'K_{\varphi}(w,z)K_{\varphi}(z,w)e^{-\varphi(w)}d\mu(w)= \nonumber \\=2Re\int K_{\varphi}(w,z)_{\varepsilon}'K_{\varphi}(z,w)e^{-\varphi(w)}d\mu(w)-\nonumber \\-\int \varphi(w)_{\varepsilon}'|K_{\varphi}(w,z)|^2e^{-\varphi(w)}d\mu(w).
\end{eqnarray}
But again by the reproducing property, $$\int K_{\varphi}(w,z)_{\varepsilon}'K_{\varphi}(z,w)e^{-\varphi(w)}d\mu(w)=K_{\varphi}(z,z)_{\varepsilon}',$$ which is real, thus the equation (\ref{Equation2}) becomes
\begin{equation}
K_{\varphi}(z,z)_{\varepsilon}'=\int \varphi(w)_{\varepsilon}'|K_{\varphi}(w,z)|^2e^{-\varphi(w)}d\mu(w).
\end{equation}

By definition we have that $$B_{\varphi}(z):=K_{\varphi}(z,z)e^{-\varphi(z)}.$$ Differentiating this with respect to $\varepsilon,$ we get
\begin{eqnarray} \label{Equation4}
B_{\varphi}(z)_{\varepsilon}':=K_{\varphi}(z,z)_{\varepsilon}'e^{-\varphi(z)}-\varphi(z)_{\varepsilon}'K_{\varphi}(z,z)e^{-\varphi(z)}= \nonumber \\=\int \varphi(w)_{\varepsilon}'|K_{\varphi}(w,z)|^2e^{-\varphi(w)-\varphi(z)}d\mu(w)- \nonumber \\-\int \varphi(z)_{\varepsilon}'|K_{\varphi}(w,z)|^2e^{-\varphi(w)-\varphi(z)}d\mu(w)=\nonumber \\=\int(\varphi(w)_{\varepsilon}'-\varphi(z)_{\varepsilon}')|K_{\varphi}(w,z)|^2e^{-\varphi(w)-\varphi(z)}d\mu(w),
\end{eqnarray}
where we have used equations (\ref{Equation3}) and (\ref{Equation2}). Since $u=\varphi_{\varepsilon}'$ we have that 
\begin{eqnarray*}
f''=(\int u B_{\varphi+\varepsilon u}d\mu)_{\varepsilon}'=\int \varphi(z)_{\varepsilon}'B_{\varphi}(z)_{\varepsilon}'d\mu (z)=\\=\int_z \int_w \varphi(z)_{\varepsilon}'(\varphi(w)_{\varepsilon}'-\varphi(z)_{\varepsilon}')|K_{\varphi}(w,z)|^2e^{-\varphi(w)-\varphi(z)}d\mu(w)d\mu (z)=\\=-\int_z \int_w (\varphi(w)_{\varepsilon}'-\varphi(z)_{\varepsilon}')^2|K_{\varphi}(w,z)|^2e^{-\varphi(w)-\varphi(z)}d\mu(w)d\mu(z)+\\+\int_z \int_w \varphi(w)_{\varepsilon}'(\varphi(w)_{\varepsilon}'-\varphi(z)_{\varepsilon}')|K_{\varphi}(w,z)|^2e^{-\varphi(w)-\varphi(z)}d\mu(w)d\mu (z)=\\=-\int_z \int_w (\varphi(w)_{\varepsilon}'-\varphi(z)_{\varepsilon}')^2|K_{\varphi}(w,z)|^2e^{-\varphi(w)-\varphi(z)}d\mu(w)d\mu(z)-\\(\textrm{by the antisymmetry in the variables } z \textrm{ and } w)\\-\int_z \int_w \varphi(z)_{\varepsilon}'(\varphi(w)_{\varepsilon}'-\varphi(z)_{\varepsilon}')|K_{\varphi}(w,z)|^2e^{-\varphi(w)-\varphi(z)}d\mu(w)d\mu (z)=\\=-\int_z \int_w (\varphi(w)_{\varepsilon}'-\varphi(z)_{\varepsilon}')^2|K_{\varphi}(w,z)|^2e^{-\varphi(w)-\varphi(z)}d\mu(w)d\mu(z)-f'',
\end{eqnarray*}
by equation (\ref{Equation4}).
This gives that $$f''=-1/2\int_z \int_w (\varphi(w)_{\varepsilon}'-\varphi(z)_{\varepsilon}')^2|K_{\varphi}(w,z)|^2e^{-\varphi(w)-\varphi(z)}d\mu(w)d\mu(z),$$ thus $f''$ is nonpositive.
\end{proof}

Here we have assumed that $\varphi(\varepsilon):=\varphi +\varepsilon u$ depends linearly on $\varepsilon.$ The general formula is 
\begin{eqnarray*}
f''=-1/2\int_z \int_w (\varphi(w)_{\varepsilon}'-\varphi(z)_{\varepsilon}')^2|K_{\varphi}(w,z)|^2e^{-\varphi(w)-\varphi(z)}d\mu(w)d\mu(z)+\\+\int \varphi(z)_{\varepsilon}''B_{\varphi}(z)d\mu(z),
\end{eqnarray*}
showing that $f$ is concave if $\varphi$ is concave in $\varepsilon.$

\begin{lemma} \label{Lemma1}
Let $f_k$ be a sequence of functions converging pointwise to a function $F$ on the interval $(-1,1).$ Furthermore, assume that for all $k,$ $f_k''\leq 0,$ and that $F'$ is continuous. Then it holds that the derivatives $f_k'$ converge pointwise to $F'.$
\end{lemma}

\begin{proof}
We pick an arbitrary point in $(-1,1),$ which we can for simplicity assume to be the origin. Since $f_k''$ is negative, $f_k'$ is decreasing. Thus for any $\varepsilon>0$ we have the inequalites $$1/\varepsilon \int_{-\varepsilon}^0 f_k' dx \geq f_k'(0) \geq 1/\varepsilon \int_0^{\varepsilon}f_k' dx,$$ i.e. $$\frac{f_k(0)-f_k(-\varepsilon)}{\varepsilon}\geq f_k'(0) \geq \frac{f_k(\varepsilon)-f_k(0)}{\varepsilon}.$$ Letting $\varepsilon$ be fix while $k$ tends to infinity we get $$\frac{F(0)-F(-\varepsilon)}{\varepsilon}\geq \limsup_{k \to \infty} f_k'(0)\geq \liminf_{k \to \infty} f_k'(0)\geq \frac{F(\varepsilon)-F(0)}{\varepsilon},$$ so by the Mean Value Theorem and the continuity of $F'$ we see that $$\lim_{k \to \infty}f_k'(0)=F'(0).$$
\end{proof}

\subsection{End of proof of Theorem 1.1}

\begin{proof}
Let us denote by $f_k$ the functions $$f_k(\varepsilon):=\frac{(n+1)!}{2k^{n+1}}\log \textrm{vol } \mathcal{B}^2(d\mu, \varphi+\varepsilon u),$$ and by $F$ the function $$F(\varepsilon):=\mathcal{E}_0(P_K(\varphi+\varepsilon u)).$$ That $f_k$ converges pointwise to $F$ is the content of Theorem 1.2 in (\cite{Berman}). By lemma (\ref{Lemma2}) we know that $f''\leq 0$. That $F'$ is continuous follows from Theorem 5.7 in (\cite{Berman}), which states that the derivative of $\mathcal{E}_0(P_K(\varphi+\varepsilon u))$ is Lipschitz. Since $$f_k'(0)=\int u \frac{B_{k\varphi}}{k^n}d\mu,$$ and $$F'(0)=\int u \textrm{MA}(P_K(\varphi)),$$ lemma (\ref{Lemma2}) and (\ref{Lemma1}) give that $$\lim_{k \to \infty} \int u \frac{B_{k\varphi}}{k^n}d\mu=\int u \textrm{MA}(P_K(\varphi)).$$ Since $u$ is an arbitrary smooth function this concludes the proof of Theorem 1.1. 
\end{proof}

E-mail addresses \\
R Berman: robertb@math.chalmers.se, D Witt Nyström: danspolitik@gmail.com

\end{document}